\documentclass[12pt]{article}
\usepackage{amsmath,amsfonts,amssymb,amsthm}

\newcommand{\I}{{\bf 1}}
\newtheorem{proposition}{Proposition}[section]
\newtheorem{theorem}[proposition]{Theorem}
\newtheorem{corollary}[proposition]{Corollary}
\newtheorem{lemma}[proposition]{Lemma}
\newtheorem{remark}[proposition]{Remark}

\newtheorem{example}[proposition]{Example}

\newcommand{\nc}{\newcommand}
\nc{\R}{{\mathbb R}}
\nc{\N}{{\mathbb N}}
\nc{\Z}{{\mathbb Z}}

\nc{\BP}{\mathbb{P}}
\nc{\BE}{\mathbb{E}}
\nc{\BQ}{\mathbb{Q}}
\nc{\bN}{{\mathbf N}}
\nc{\BX}{{\mathbb X}}
\nc{\BY}{{\mathbb Y}}
\nc{\bM}{{\mathbf M}}
\nc{\bF}{{\mathbf F}}
\nc{\bG}{{\mathbf G}}
\nc{\bH}{{\mathbf H}}
\nc{\bD}{{\mathbf D}}

\begin{document}

\author{G\"unter Last\footnote{
Institut f\"ur Stochastik,  Karlsruher Institut f\"ur Technologie,
76128 Karlsruhe, Germany. 
Email: guenter.last@kit.edu}
\ and Mathew D.\ Penrose 
\footnote{
Department of Mathematical Sciences, University of Bath,
Bath BA2 7AY, United Kingdom,
Email: m.d.penrose@bath.ac.uk} 
\footnote{Partially supported by
the Alexander von Humboldt Foundation through
a Friedrich Wilhelm Bessel Research Award.}
}
\title{Martingale representation for Poisson processes\\
with applications to minimal variance hedging}
\date{\today}
\maketitle
\begin{abstract}
\noindent

We consider a Poisson process $\eta$ on a 
measurable space $(\BY,\mathcal{Y})$ equipped with
a partial ordering, assumed to be strict almost everwhwere
with respect to the intensity measure $\lambda$ of $\eta$.
We give a Clark-Ocone type formula providing an explicit  representation 
of square integrable martingales (defined with respect to the
natural filtration associated with $\eta$), which was previously 
known only in the special case, when $\lambda$
is the product of Lebesgue measure on $\R_+$ and
a $\sigma$-finite measure on another space $\BX$. Our proof
is new and based on only a few basic properties of 
Poisson processes and stochastic integrals.
We also consider the more general case of an
independent random measure in the sense of It\^o of pure jump type
and show that the Clark-Ocone type representation leads to an
explicit version of the Kunita-Watanabe decomposition
of square integrable martingales.
We also find the explicit minimal variance
hedge in a quite general financial market driven by an independent random measure.
\end{abstract}


\vspace*{0.1cm}

\section{Introduction}
\setcounter{equation}{0}

Any square integrable martingale with respect to a Brownian
filtration can be written as a stochastic integral,
see \cite{Ito51} and Theorem 18.10 in \cite{Kallenberg}.
This {\em martingale representation theorem}
is an important result of stochastic analysis.
Similar results are available for marked point processes
(see e.g.\ \cite{LB95,JaSh03} and the references given there)
and for general semimartingales, see Section III.4 in 
\cite{JaSh03}.
For some Brownian martingales Clark \cite{Clark70} found a more
explicit version of the integrand in the representation.
Ocone \cite{Ocone84} revealed the relationship of Clark's formula
to Malliavin calculus.

The topic of the present paper is a Clark-Ocone type martingale
representation formula when the underlying 
filtration is generated by a Poisson process $\eta$ on
a measurable space $(\BY,\mathcal{Y})$ equipped with
a partial ordering.
Our main result (Theorem \ref{tmrt}) provides
a representation of square integrable martingales
as a (stochastic) Kabanov-Skorohod integral with respect to the
compensated Poisson process. 
In the case $\BY=\R_+\times\BX$ is the product of
$\R_+:=[0,\infty)$ and a Borel space $\BX$,
special cases of this formula are well-known.
Stationary  Poisson processes on $\R_+$ were treated in 
Picard \cite{Pic96}, while \cite{AneLed00} considered
the more general case of a finite set $\BX$.
In \cite{Wu00} it was shown how to use the
Malliavin  calculus for Poisson processes developed 
in \cite{Ogura72,Kab75,NuViv90} and the results in \cite{DKW88}
to get the Clark-Ocone formula under an
additional integrability assumption 
in the case where the intensity measure of $\eta$ 
is the product of Lebesgue measure and a $\sigma$-finite
measure on $\BX$. This is also the approach taken in 
\cite{Lo05} and \cite{DiNOkPro09}
when treating  pure jump L\'evy processes 
(without refering to \cite{Wu00}).
Translated to our setting, this is again the special case where 
the intensity measure has product form. 
Our proof of Theorem \ref{tmrt} is based on
the explicit Fock space representation of Poisson functionals
\cite[Theorem 1.5]{LaPe09} and the basic isometry properties of
stochastic integrals,
and is distinct from the proofs of related results
that we have seen in the literature.
In particular we are not using any other
martingale representation theorem for Poisson spaces.

We apply Theorem \ref{tmrt} to 
derive the explicit Kunita-Watanabe projection
of square integrable martingales onto the space
of stochastic integrals against an independent
random measure (in the sense of It\^o \cite{Ito56}) without Gaussian
component. We also find the explicit minimal variance
hedge in a quite general market driven by an independent
random measure.

We now describe the contents of this paper in more detail.
Throughout the paper we consider a Poisson process $\eta$
on a measurable space $(\BY,\mathcal{Y})$
with $\sigma$-finite intensity measure $\lambda$.
The underlying probability space is denoted by $(\Omega,\mathcal{F},\BP)$.
We can interpret $\eta$ as a random element in the space
$\bN:=\bN(\BY)$ of $\sigma$-finite integer-valued measures $\mu$
on $\BY$ equipped with the smallest $\sigma$-field
making the mappings $\mu\mapsto\mu(B)$ measurable for all
$B\in\mathcal{Y}$. We assume that $\BY$ is equipped with
a transitive binary relation $<$ such that
$\{(y,z):y<z\}$ is a measurable subset of $\BY^2$
and such that for any $y,z\in\BY$ at most one of the
relations $y<z$ and $z<y$ can be satisfied. 
We also assume
that $<$ strictly orders the points of $\BY$ $\lambda$-a.e., that is
\begin{align}\label{diffuse}
\lambda([y])=0, \quad y\in\BY,
\end{align}
where $[y]:=\BY\setminus \{z\in\BY:\text{$z<y$ or $y<z$}\}$.
For any $\mu\in\bN$ let $\mu_y$ denote the restriction
of $\mu$ to 
$y_{\downarrow}:=\{z\in\BY:z<y\}$. Our final assumption on
$<$ is that $(\mu,y)\mapsto \mu_y$ is a measurable mapping
from $\bN\times\BY$ to $\bN$.

For $y\in\BY$ the {\em difference operator} $D_y$
is given as follows.
For any measurable $f:\bN\rightarrow\R$ the function
$D_yf$ on $\bN$ is defined by
\begin{align}\label{addone}
D_{y}f(\mu):=f(\mu+\delta_{y})-f(\mu),\quad \mu\in\bN,
\end{align}
where $\delta_{y}$ is the Dirac measure located at a point $y\in\BY$.
We need a version of the
conditional expectation $\BE[D_yf(\eta)|\eta_y]$
that is jointly measurable in all arguments.
Thanks to the independence properties of a Poisson process
we can and will work with
\begin{align}\label{956}
\BE[D_yf(\eta)|\eta_y]:=\int D_yf(\eta_y+\mu)\Pi^y(d\mu),
\end{align}
where $\Pi^y$ is the distribution of the restriction
of $\eta$ to $\BY\setminus y_{\downarrow}$. We use this
definition only if the right-hand side is well defined and finite.
Otherwise we set $\BE[D_yf(\eta)|\eta_y]:=0$. Note that
$\BE[D_yf(\eta)|\eta_y]=h(\eta,y)$,
where $h:\bN\times\BY\rightarrow\R$ is defined by
\begin{align}\label{h}
h(\mu,y):=\int D_yf(\mu_y+\nu)\Pi^y(d\nu). 
\end{align}
Since $(\mu,y)\mapsto \mu_y$ is assumed measurable,
the function $h$ is measurable as well. Moreover, it satisfies
\begin{align}\label{predict}
h(\mu,y)=h(\mu_y,y),\quad (\mu,y)\in\bN\times\BY.
\end{align}
Justified by Proposition \ref{l1}
we call a measurable function $h$ with the property 
\eqref{predict} {\em predictable}, see Remark \ref{rmnot}.
This notion depends on the ordering $<$. The fact that
this dependence is not reflected in our terminology,
will not lead to confusion.

If $h:\bN\times\BY\rightarrow\R$ is a measurable function
then we denote by $\delta(h)\equiv\int h(\eta,y)\hat\eta(dy)$ the stochastic 
Kabanov-Skorohod integral of $h$ with respect to the compensated
Poisson process $\hat\eta:=\eta-\lambda$ \cite{Kab75,Skor75,KabSk75}. 
This integral is well defined
only, if the integrability condition \eqref{domainS} on $h$ holds.
If, in addition, 
$h\in L^1(\BP_\eta\otimes\lambda)\cap L^2(\BP_\eta\otimes\lambda)$,
then Theorem 3.5 in \cite{LaPe09} provides a pathwise
interpretation of $\delta(h)$:
\begin{align}\label{tequal}
\delta(h)=\int h(\eta -\delta_y,y)\eta(dy)
-\int h(\eta,y)\lambda(dy) \quad \BP\text{-a.s.}
\end{align}
In fact, if $h\in L^2(\BP_\eta\otimes\lambda)$ is predictable
(i.e.\ \eqref{predict} holds), then $\delta(h)$ is well defined
and we have the isometry relation 
\begin{align}\label{l2i3}
\BE \delta(h)^2=\BE\int h(\eta,y)^2\lambda(dy).
\end{align}
We prove these facts in Section \ref{srepP}, see 
Propositions \ref{text} and \ref{pquadcha}. For
predictable functions $h$
equation \eqref{l2i3} can be used to extend \eqref{tequal}
from $L^1(\BP_\eta\otimes\lambda)\cap L^2(\BP_\eta\otimes\lambda)$
to $L^2(\BP_\eta\otimes\lambda)$.
If $h\in L^2(\BP_\eta\otimes\lambda)$ is predictable
and $A\in\mathcal{Y}$, then we can define 
$\int_A h(\eta,y)\hat\eta(dy):=\delta(\I_{\bN\times A} h)$.
Let $\BP_\eta$ denote the distribution of $\eta$.
For $f\in L^2(\BP_\eta)$ (i.e.\ for measurable $f:\bN\rightarrow\R$
with $\BE f(\eta)^2<\infty$) we have the following representation
of $f(\eta)$.

\begin{theorem}\label{tmrt} Let $\eta$ be a Poisson process
on $\BY$ with an intensity measure $\lambda$
satisfying \eqref{diffuse} and let $f\in L^2(\BP_\eta)$.
Then
\begin{align}\label{1}
\BE\int\BE[D_yf(\eta)|\eta_y]^2\lambda(dy)<\infty
\end{align}
and we have $\BP$-a.s.\ that 
\begin{align}\label{mr}
f(\eta)=\BE f(\eta)+\int \BE[D_yf(\eta)|\eta_y]\hat\eta(dy).
\end{align}
Moreover, we have for any $y\in\BY$ that $\BP$-a.s.
\begin{align}\label{mr47}
\BE[f(\eta)|\eta_y]=\BE f(\eta)+
\int_{y_{\downarrow}}\BE[D_zf(\eta)|\eta_z]\hat\eta(dz).
\end{align}
\end{theorem}

Define $M_y:=\BE[f(\eta)|\eta_y]$, $y\in\BY$, where $f$
is as in \eqref{mr47}. If $z<y$ then the $\sigma$-field
$\sigma(\eta_z)$ is contained in $\sigma(\eta_y)$ and we have
the martingale property $\BE[M_y|\eta_z]=M_z$ a.s.
Equation \eqref{mr47} provides an explicit representation
of the martingale $(M_y)$ as stochastic integral
of an explicitly known integrand.

In the remainder of this introduction we
assume that $\BY=\R_+\times\BX$, where $(\BX,\mathcal{X})$ is
a Borel space and that $(s,x)<(s',x')$ if and only if $s<s'$. 
Assumption \eqref{diffuse} means that 
\begin{align}\label{diffuse2}
\lambda(\{t\}\times\BX)=0,\quad t\ge 0.
\end{align}
We do not assume $\lambda$ to be of product form.
In Section \ref{spred} we first discuss Theorem \ref{tmrt}
in this case. Then we show that a function is predictable
essentially if and only if it is predictable in the
standard sense of stochastic analysis.
Theorem \ref{text2} shows that the Kabanov-Skorohod integral
of a predictable function coincides with the 
(standard) stochastic integral. This extends
results in \cite{Kab75} and \cite{NuViv90} for Poisson
processes on $\R_+$.

In Section \ref{sl2} we consider instead of the compensated
Poisson process $\hat\eta$ a more general centred independent
random measure $\zeta$ (in the sense of \cite{Ito56})
on $\R_+\times\BX$. We assume that $\zeta$ has no Gaussian part and a
$\sigma$-finite variance measure with diffuse projection
onto the first coordinate. Then $\zeta$ can be
represented in terms of a Poisson process $\eta$ as above
on $\BY:=\R_+\times\BX\times(\R\setminus\{0\})$.
Consequently we can apply our Clark-Ocone type formula to
obtain an explicit formula for the orthogonal projection
of a square integrable function of $\eta$ onto the space
of all stochastic integrals against $\zeta$, see Theorem \ref{tmrtKW}.
Such projections were first considered
by Kunita and Watanabe \cite{KuWa67} in the setting
of continuous martingales. Later these ideas were extended
to semimartingales, see e.g.\ Schweizer \cite{schweizer94}. 
Using a different approach (and allowing for a Gaussian component)
Di Nunno  \cite{Nunno07} proved a version of Theorem \ref{tmrtKW} for 
special (``core'') functions of $\eta$. 
In fact we prove our results in the more general case
of an independent random measure $\zeta$ on a Borel space
$(\BY',\mathcal{Y}')$ with a diffuse
and $\sigma$-finite variance measure $\beta$ such that 
$\BY'$ is ordered almost everwhere with respect to $\beta$.

In Section \ref{smeanva} we consider
a quite general financial market with a continuum of assets,
driven by an independent random measure without Gaussian
component. Again all processes can be
represented in terms of a Poisson process $\eta$
on a suitable state space. A function $f\in L^2(\BP_\eta)$
can then be interpreted as a contingent claim.
Minimizing the $L^2$-distance 
between $f(\eta)-\BE f(\eta)$ and a certain space of
stochastic integrals against the assets, yields the
{\em minimal variance hedge} of $f(\eta)$. 
Theorem \ref{tminvar} finds this hedge explicitly,
while Theorem \ref{texacth} identifies the claims that
can be perfectly hedged. These theorems extend the main results 
in \cite{BNLOk03}, which treats the case of a market
driven by a finite number of independent L\'evy processes.


\section{Representation of Poisson martingales}\label{srepP}
\setcounter{equation}{0}

In this section we prove Theorem \ref{tmrt}, starting with some
definitions and preliminary observations. 
Let $f:\bN\rightarrow\R$ be a measurable function.
For $n\ge 2$ and $(y_1,\ldots,y_n)\in\BY^n$
we define a function
$D^{n}_{y_1,\ldots,y_n}f:\bN\rightarrow\R$
inductively by
\begin{align}\label{differn}
D^{n}_{y_1,\ldots,y_{n}}f:=D^1_{y_{1}}D^{n-1}_{y_2,\ldots,y_{n}}f,
\end{align}
where $D^1:=D$ and $D^0f=f$. For $f\in L^2(\BP_\eta)$
it was proved in \cite{LaPe09}
that $D^{n}_{y_1,\ldots,y_{n}}f(\eta)$ is integrable for
$\lambda^n$-a.e.\ $(y_1,\ldots,y_n)$ and that
\begin{align}\label{Tn}
T_n f(y_1,\ldots,y_n):=\BE D^n_{y_1,\ldots,y_n} f(\eta),
\quad (y_1,\ldots,y_n)\in\BY^n,
\end{align}
defines a symmetric function in $L^2(\lambda^n)$.
Moreover, we have the Wiener-It\^o chaos expansion
\begin{align}\label{chaos}
f(\eta)=\sum^\infty_{n=0}\frac{1}{n!}I_n(T_nf),
\end{align}
where the series converges in $L^2(\BP)$. 
Here $I_n(g)$ denotes the $n$th {\em multiple Wiener-It\^o integral}
of a symmetric $g\in L^2(\lambda^n)$, see \cite{Ito56}.
These integrals satisfy the orthogonality relations
\begin{align}\label{orth}
\BE I_m(g)I_n(h)=\I\{m=n\}m!\langle g , h \rangle_n
\quad m,n\in\N_0,
\end{align}
where $\langle\cdot,\cdot\rangle_n$ denotes the scalar product
in $L^2(\lambda^n)$.

Let $h\in L^2(\BP_\eta\otimes\lambda)$. Then $h(\cdot,y)\in L^2(\BP_\eta)$
for $\lambda$-a.e.\ $y$ and we may consider the chaos expansion 
\begin{align}\label{ch}
h(\eta,y)=\sum^\infty_{n=0} I_n(h_n(y)),
\end{align}
where $h_n(y)\in L^2(\lambda^n)$, $n\in\N$, are given by
\begin{align}\label{ch9}
h_n(y)(y_1,\ldots,y_n):=\BE D^n_{y_1,\ldots,y_n} f(\eta,y).
\end{align}
Let $\tilde h_n$ be the {\em symmetrization} of
this function, that is
\begin{align*}
\tilde h_n(y_1,\ldots,y_{n+1})=
\frac{1}{n+1}
\sum^n_{i=1}h_n(y_{i})(y_1,\ldots,y_{i-1},y_{i+1},\ldots,y_{n+1}).
\end{align*}
From \eqref{ch} and \eqref{orth}  we obtain
that $\tilde h_n\in L^2(\lambda^{n+1})$ and we can define
the {\em Kabanov-Skorohod integral} \cite{Hitsuda72,Kab75,Skor75,KabSk75,LaPe09}
of $h$, denoted $\delta(h)$, by
\begin{align}\label{Skorint}
\delta(h):=\sum^\infty_{n=0} I_{n+1}(\tilde h_n),
\end{align}
which converges in $L^2(\BP)$ provided that
\begin{align}\label{domainS}
\sum^\infty_{n=0}(n+1)!\int \tilde h_n^2d\lambda^{n+1}<\infty.
\end{align}
We need the following duality relation
from \cite{NuViv90}, see also Proposition 3.4 in \cite{LaPe09}. 
We let $\|\cdot\|_n$ denote the norm in $L^2(\lambda^n)$. 

\begin{proposition}\label{Spath} Assume that $g\in L^2(\BP_\eta)$
satisfies 
\begin{align}\label{conv}
\sum^\infty_{n=1}\frac{1}{(n-1)!}\|T_ng\|^2_n<\infty.
\end{align}
Let $h\in L^2(\BP_\eta\otimes\lambda)$
with a chaos expansion satisfying \eqref{domainS}. Then
$\BE\int (D_yg(\eta))^2\lambda(dy)<\infty$ and
\begin{align}\label{adj2}
\BE\int D_yg(\eta)h(\eta,y)\lambda(dy)=\BE g(\eta)\delta(h).
\end{align}
\end{proposition}

Proposition \ref{Spath} easily shows that $\delta$ is
{\em closed}, see \cite{Kab75} and \cite{NuViv90}.
This means that if $h_k\in L^2(\BP_\eta\otimes\lambda)$, $k\in\N$,
satisfy \eqref{domainS}, $h_k\to h$ in $L^2(\BP_\eta\otimes\lambda)$
and $\delta(h_k)\to X$ in $L^2(\BP)$, then $h$
satisfies \eqref{domainS} and $\delta(h)=X$ a.s.
We shall use this fact repeatedly in the sequel.

The next result shows that the Kabanov-Skorohod integral of
a predictable $h$ is defined, if $h$
is square integrable with respect to $\BP_\eta\otimes\lambda$.

\begin{proposition}\label{text} 
Let $h\in L^2(\BP_\eta\otimes\lambda)$ be predictable.
Then \eqref{domainS} holds.
\end{proposition}
{\sc Proof:} Consider the functions defined by \eqref{ch9}.
Since $h$ is predictable, we have that $h_n(y)(y_1,\ldots,y_n)=0$
whenever $y_i>y$ for some $i\in\{1,\ldots,n\}$.
This implies that 
$$
\I_{\Delta_{n+1}}(y_1,\ldots,y_{n+1})\tilde h_n(y_1,\ldots,y_{n+1})
=\I_{\Delta_{n}}(y_1,\ldots,y_n)\frac{1}{n+1}
h_n(y_{n+1})(y_1,\ldots,y_n),
$$
where 
\begin{align}\label{89}
\Delta_n:=\{(y_1,\ldots,y_n)\in \BY^n:y_1<\ldots <y_n\}.
\end{align}
In view of \eqref{diffuse} it follows that
\begin{align*}
\|\tilde h_n\|^2_{n+1}&=(n+1)!\|\I_{\Delta_{n+1}}\tilde h_n\|^2_{n+1}\\
&=\frac{(n+1)!}{(n+1)^2}\int\|\I_{\Delta_n}h_n(y)\|^2_n\lambda(dy)
=\frac{1}{n+1}\int\|h_n(y)\|^2_n\lambda(dy).
\end{align*}
Hence we obtain from \eqref{orth} and \eqref{ch} that
 \begin{align*}
\sum^\infty_{n=0}(n+1)!\|\tilde h_n\|^2_{n+1}
&=\sum^\infty_{n=0}\int n!\|h_n(y)\|^2_n\lambda(dy)\\
&=\sum^\infty_{n=0}\int \BE I_n(h_n(y))^2\lambda(dy)
=\int\BE h(y)^2\lambda(dy)<\infty.
\end{align*}
Therefore \eqref{domainS} holds.\qed

\vspace{0.3cm}
Let $h\in L^2(\BP_\eta\otimes\lambda)$ be predictable
and $B\in\mathcal{Y}$. Then
$\I_{\bN\times B}h\in L^2(\BP_\eta\otimes\lambda)$
is also predictable. Moreover, we have from  \eqref{adj2} that
\begin{align}\label{zero}
\delta(\I_{\bN\times B}h)=0\quad \BP\text{-a.s. if $\lambda(B)=0$}.
\end{align}
The following proposition implies a part of Theorem \ref{tmrt}.

\begin{proposition}\label{pmart}
Let $h\in L^2(\BP_\eta\otimes\lambda)$ be predictable.
Then, for any $y\in \BY$,
\begin{align}\label{My}
\BE\Big[\int h(\eta,z)\hat\eta(dz)\Big|\eta_y\Big]=
\int_{y_{\downarrow}}h(\eta,z)\hat\eta(dz) \quad \BP\text{-a.s.}
\end{align}
\end{proposition}
{\sc Proof:} The right-hand side of \eqref{My}
can be chosen $\sigma(\eta_y)$-measurable. This fact can be
traced back to \eqref{orth}: if $f\in L^2(\lambda^n)$ is symmetric
and vanishes outside $B^n$ for some $B\in\mathcal{Y}$
and $I_n^B$ denotes the $n$th Wiener-Ito integral
with respect to the restriction of $\eta$ to $B$,
then $I_n(f)=I^B_n(f)$ $\BP$-a.s. 

To prove \eqref{My}, we take $y\in\BY$ and a measurable function
$g:\bN\rightarrow\R$ such that the function $g_y$ defined
by $g_y(\mu):=g(\mu_y)$ satisfies \eqref{conv}.
Since $D_zg_y=0$ for $y<z$ we obtain from Proposition \ref{Spath}
that
$$
0=\BE g(\eta_y)\int\I\{y<z\}h(\eta,z)\hat\eta(dz).
$$
From \eqref{zero} and \eqref{diffuse} we have
\begin{align}\label{87}
\int\I_{[y]}(z)h(\eta,z)\hat\eta(dz)=0\quad \BP\text{-a.s.}
\end{align}
Hence we obtain from the linearity of $\delta$ that
\begin{align}\label{2.16}
\BE g(\eta_y)\int h(\eta,z)\hat\eta(dz) 
=\BE g(\eta_y)\int_{y_{\downarrow}}h(\eta,z)\hat\eta(dz).
\end{align}
Now we consider a function $g$ of the
form $g(\mu):=\exp[-\int hd\mu]$, where $h:\BY\rightarrow\R_+$
is measurable and vanishes outside a set $C\in\mathcal{Y}$
with $\lambda(C)<\infty$. It can be easily checked, that
$g_y$ satisfies \eqref{conv} (cf.\ also the proof of Theorem 3.3
in \cite{LaPe09}). Hence \eqref{2.16} holds
for all linear combinations of such functions.
A monotone class argument shows that \eqref{2.16}
holds for all bounded measurable $g:\bN\rightarrow\R$
(cf.\ the proof of Lemma 2.2 in \cite{LaPe09}).
This is enough to deduce \eqref{My}.\qed

\vspace{0.3cm}
\noindent
{\sc Proof of Theorem \ref{tmrt}:} Let $f\in L^2(\BP_\eta)$
and define $h:\bN\times\BY\rightarrow\R$ by \eqref{h}.
Then $h$ is predictable. 
Moreover, Theorem 1.5 in \cite{LaPe09}
implies that $h\in L^2(\BP_\eta\otimes\lambda)$, that is
\eqref{1} holds. By Proposition \ref{text}, the Kabanov-Skorohod integral
$\delta(h)$ is well defined. We have to show that 
\begin{align}\label{3789}
f(\eta)=\BE f(\eta)+\delta(h)\quad \BP\text{-a.s.}
\end{align}
Let $g\in L^2(\BP_\eta)$ satisfy \eqref{conv}.
By Proposition \ref{Spath}, 
\begin{align*}
\BE g(\eta)\delta(h)
&=\BE\int D_yg(\eta)\BE[D_yf(\eta)|\eta_y]\lambda(dy)\\
&=\int\BE[\BE[D_yg(\eta)|\eta_y]\BE[D_yf(\eta)|\eta_y]]
\lambda(dy),
\end{align*}
where the second equality comes from Fubini's theorem
and a standard property of conditional expectations.
Applying Theorem 1.5 in \cite{LaPe09}, we obtain that
\begin{align*}
\BE g(\eta)\delta(h)=\BE g(\eta)f(\eta)-(\BE g(\eta))(\BE f(\eta)),
\end{align*}
that is $\BE g(\eta)(\BE f(\eta)+\delta(h))=\BE g(\eta)f(\eta)$.
Since the set of all $g\in L^2(\BP_\eta)$ satisfying \eqref{conv}
is dense in $L^2(\BP_\eta)$, we obtain \eqref{3789}.
The remaining assertion follows from Proposition \ref{pmart}. \qed

\vspace{0.3cm}
We finish this section with a standard property
of stochastic integrals.

\begin{proposition}\label{pquadcha} 
Let $h,\tilde{h}\in L^2(\BP_\eta\otimes\lambda)$ be predictable.
Then 
\begin{align}\label{l2i}
\BE \delta(h)\delta(\tilde{h})=
\BE\int h(\eta,y)\tilde{h}(\eta,y)\lambda(dy).
\end{align}
\end{proposition}
{\sc Proof:} By linearity and polarization we may
assume that $h=\tilde{h}$. Let us first assume that
$h$ is bounded and that $h(\mu,x)=0$ for 
$x\notin C\in\mathcal{Y}$, where $\lambda(C)<\infty$. 
In particular,
$h\in L^p(\BP_\eta\otimes\lambda)$ for any $p>0$.
By \eqref{tequal}, 
\begin{align}\label{890} \notag
\BE \delta(h)^2=&\BE\left(\int h(\eta-\delta_y,y)\eta(dy)\right)^2\\
&-2\BE\left(\int h(\eta-\delta_x,x)\eta(dx)\int h(\eta,y)\lambda(dy)\right)
+\BE\left(\int h(\eta,y)\lambda(dy)\right)^2.
\end{align}
Our assumptions on $h$ guarantee that all these
expectations are finite.
We are now performing a fairly standard calculation
based on the Mecke equation, see e.g.\ (2.10) in \cite{LaPe09}.
The first term on the right-hand side of \eqref{890}
equals
\begin{align*}
\BE\int &h(\eta,y)^2\lambda(dy)
+\BE\iint h(\eta+\delta_y,x)h(\eta+\delta_x,y)\lambda(dy)\lambda(dx)\\
=&\BE\int h(\eta,y)^2\lambda(dy)
+2\BE\iint \I\{x<y\}h(\eta,x)h(\eta+\delta_x,y)\lambda(dy)\lambda(dx),
\end{align*}
where we have used symmetry, \eqref{diffuse} and \eqref{predict},
to obtain the equality.
The second term on the right-hand side of \eqref{890} equals
\begin{align*}
-2\BE\iint \I\{x<y\} &h(\eta,x)h(\eta+\delta_x,y)\lambda(dy)\lambda(dx)\\
&-2\BE\iint \I\{y<x\} h(\eta,x)h(\eta,y)\lambda(dy)\lambda(dx).
\end{align*}
Summarizing, we obtain that \eqref{l2i3} holds, as required.

In the general case we define, for $k\in\N$,
$$
h_k(\mu,x):=\I\{|h(\mu,x)|\le k\}\I\{x\in C_k\}h(\mu,x),
\quad (\mu,x)\in\bN\times\BY,
$$
where $C_k\uparrow\BY$ and $\lambda(C_k)<\infty$.
The functions $h_k$ are predictable and satisfy the
assumptions made above. From dominated
convergence we have 
$\BE\int(h(\eta,x)-h_k(\eta,x))^2\lambda(dx)\to 0$
as $k\to\infty$.
Then \eqref{l2i3} implies that $\delta(h_k)$
is a Cauchy sequence in $L^2(\BP)$ and hence converges
towards some $X\in L^2(\BP)$. Since $\delta$ is closed,
we obtain $X=\delta(h)$ and hence the assertion.\qed

\section{Martingales and stochastic integration}\label{spred}
\setcounter{equation}{0}

Assume that $\BY=\R_+\times\BX$, where
$(\BX,\mathcal{X})$ is a measurable space. We define
$(s,x)<(s',x')$ if and only if $s<s'$. 
Throughout this section we consider a Poisson process
$\eta$ on $\BY$ whose intensity measure $\lambda$ is 
$\sigma$-finite and satisfies \eqref{diffuse2}.
We discuss Theorem \ref{tmrt}
and the Kabanov-Skorohod integral of predictable functions.

For any $s\ge 0$ and $\mu\in\bN$ we denote by
$\mu_s$ (resp.\ $\mu_{s-}$) the restriction of $\mu$ to 
$[0,s]\times\BX$ (resp.\ $[0,s)\times\BX$).
Theorem \ref{tmrt} takes the following form.

\begin{theorem}\label{tmrt2} 
Let $f\in L^2(\BP_\eta)$. Then
\begin{align}\label{125}
\BE\int\BE[D_{(s,x)}f(\eta)|\eta_{s-}]^2\lambda(d(s,x))<\infty
\end{align}
and we have for any $t\ge 0$ that $\BP$-a.s.\ 
\begin{align}\label{mr34}
\BE[f(\eta)|\eta_t]=
\BE f(\eta)+
\int \I_{[0,t]}(s) \BE[D_{(s,x)}f(\eta)|\eta_{s-}]\hat\eta(d(s,x)).
\end{align}
\end{theorem}
{\sc Proof:} Relation \eqref{125} follows directly from
Theorem \ref{tmrt}.
For any $t\ge 0$ we have $\BP$-a.s.\ that
$$
\BE[f(\eta)|\eta_t]=\int f(\eta_{t}+\mu)\Pi^t(d\mu)
$$
where $\Pi^t$ is the distribution of the restriction
of $\eta$ to $(t,\infty)\times\BX$, compare with \eqref{956}.
By \eqref{diffuse2}, $\Pi^t$ is also the distribution of the restriction
of $\eta$ to $[t,\infty)\times\BX$ and $\eta_t=\eta_{t-}$ a.s.
Hence $\BE[f(\eta)|\eta_t]=\BE[f(\eta)|\eta_{t-}]$ and 
\eqref{mr34} follows from \eqref{mr47} and \eqref{87}.\qed

\begin{remark}\label{rpmart}\rm 
Let $h\in L^2(\BP_\eta\otimes\lambda)$ be predictable
and define 
$$
M_t:=\int\I_{[0,t]}(s)h(\eta,s,x)\hat\eta(d(s,x)),\quad
t\in[0,\infty].
$$
Proposition \ref{pmart} and \eqref{87} imply for any
$t\in[0,\infty]$ that $\BE[M_\infty|\eta_{t-}]=M_t$ $\BP$-a.s. 
In the proof of Theorem \ref{tmrt2} we have seen that
$\BE[M_\infty|\eta_{t-}]=\BE[M_\infty|\eta_t]$ $\BP$-a.s.
Hence $(M_t)_{t\in[0,\infty]}$ is a martingale with respect to the filtration
$(\sigma(\eta_t))_{t\in[0,\infty]}$, where $\eta_\infty:=\eta$.
This martingale is {\em square integrable}, that is
$M_\infty\in L^2(\BP)$.
\end{remark}

Our next aim is to clarify the meaning of
the predictability property \eqref{predict} and
to discuss the Kabanov-Skorohod integral of predictable functions.
To do so, we introduce a measurable subset
$\bN^*$ of $\bN$ as follows. Let $C_1,C_2,\ldots$ be a sequence
of disjoint measurable subsets of $\BY$ with union $\BY$.
We let $\bN^*$ be the set of all $\mu\in\bN$ having the properties
$\mu(\{0\}\times\BX)=0$ and
$\mu(C_n)<\infty$ for all $n\in\N$. 
For any $t\in[0,\infty]$ let $\mathcal{N}_t$ the smallest
$\sigma$-field of subsets of $\bN^*$, making the mappings 
$\mu\mapsto\mu(B\cap ([0,t]\times\BX))$ measurable for all
$B\in\mathcal{Y}$. Here $\mu_\infty:=\mu$.
The {\em predictable $\sigma$-field}  $\mathcal{P}$ 
(see \cite{JaSh03}) is the
smallest $\sigma$-field containing the sets
\begin{align}\label{pred2}
A\times(s,t]\times B,\quad s<t,A\in \mathcal{N}_s,B\in\mathcal{X}.
\end{align}

The next proposition provides a useful characterization
of the predictable $\sigma$-field. We have to assume that
$(\BX,\mathcal{X})$ is Borel isomorphic to a Borel subset of $[0,1]$.
Such a space is called {\em Borel space}, see \cite{Kallenberg}.

\begin{proposition}\label{l1} Assume that $(\BX,\mathcal{X})$ is
a Borel space. Let
$h:\bN^*\times\R_+\times\BX\rightarrow\R$ be
measurable. Then $h$ is $\mathcal{P}$-measurable 
if and only if \eqref{predict} holds, that is
\begin{align}\label{predict2}
h(\mu,s,x)=h(\mu_{s-},s,x),\quad (\mu,s,x)\in\bN^*\times\BX\times\R_+.
\end{align}
\end{proposition}
{\sc Proof:} 
The filtration $(\mathcal{N}_t)_{t\ge 0}$
is not right-continuous, but has otherwise many
of the properties of a point process filtration as
studied in Section 2.2 of \cite{LB95}. 
To make this more precise, we introduce
$\bN_n$, $n\in\N$, as the set of all finite integer-valued measures $\mu$
on $C_n$ such that $\mu((\{0\}\times\BX)\cap C_n)=0$.
Any $\mu\in\bN_n$ can be written as
\begin{align}\label{sxmu}
\mu(B)=\sum^m_{i=1}\int\I_B(s_i,x)\mu_i(dx),
\end{align}
where $m\ge 0$, $0<s_1<\ldots<s_m$, and $\mu_1,\ldots,\mu_m$,
are finite non-trivial integer-valued measures on $\BX$. 
(Here we use the Borel structure of $\BX$.)
It is convenient to identify $\mu$ with the infinite sequence
$(s_i,\mu_i)$, $i\in\N$, where $(s_i,\mu_i):=(\infty,0)$ for $i>m$, and
$0$ denotes the zero measure. Let $\bN_{f}(\BX)$ be the space
of all finite counting measures on $\BX$. It is not difficult to
see the this is a Borel space. Moreover, 
the quantities $m$, $s_1,\ldots,s_m$, $\mu_1,\ldots,\mu_m$ in \eqref{sxmu}
depend on $\mu$ in a measurable way.
(This does require the Borel structure of $\BX$ and $\bN_f(\BX)$). 
Therefore we can identify
$\bN_n$ with a measurable subset $\bN'_n$ of the space $\bM$ defined as 
the set of all sequences 
$((s_i,\mu_i))_{i\in\N}\in ((0,\infty]\times \bN_f(\BX))^\infty$
with the following properties. If $s_i<\infty$, then $s_i<s_{i+1}$
and $\mu_i\ne 0$. If  $s_i=\infty$, then $s_{i+1}=\infty$ and $\mu_i=0$.
The space $\bN'_n\subset\bM$ can be equipped with the product
topology inherited from $([0,\infty]\times \bN_f(\BX))^\infty$.
Now we indentify the whole space $\bN^*$ with 
$\bN'_1\times\bN'_2\times\ldots$, again equipped with the product topology.
The crucial property
of this topology is that the mappings $s\mapsto \mu_s$ and
$s\mapsto \mu_{s-}$ are right-continuous respectively left-continuous.
Therefore it is not difficult to check that
Theorem 2.2.6 in \cite{LB95} applies to the filtration $(\mathcal{N}_t)$.
\qed

\begin{remark}\label{rmnull}\rm The assumption $\mu(\{0\}\times\BX)=0$
for $\mu\in \bN^*$ has been made for convenience.
Without this condition the $\sigma$-field $\mathcal{N}_0$
becomes non-trivial, and we have to include the sets
$A\times\{0\}\times B$ ($A\in \mathcal{N}_0$, $B\in\mathcal{X}$)
into the $\sigma$-field $\mathcal{P}$. 
If we then redefine $\mu_{0-}$ as the
restriction of $\mu$ to $\{0\}\times\BX$, Proposition \ref{l1}
remains valid.
\end{remark}

\vspace*{0.3cm}
We now assume that the sets $C_n$, $n\in\N$,
are chosen in such a way, that the
intensity measure $\lambda$ of $\eta$ is 
finite on these sets. Let $\eta^*$ be the random
element in $\bN^*$, defined by $\eta^*:=\eta$ if
$\eta\in N^*$ and $\eta^*:=0$, otherwise. The second case
has probability 0. Let $F^*_1$ and $F^*_2$ denote the
$\mathcal{P}$-measurable elements of 
$L^1(\BP_{\eta^*}\otimes\lambda)$ and $L^2(\BP_{\eta^*}\otimes\lambda)$
respectively.
For $h\in F^*_2$ we can define the stochastic integral
$\delta^*(h)$ of $h$ against the compensated Poisson process
$\eta^*-\lambda$ in the following standard way, see e.g.\ \cite{IkWa81}.
If $h\in F^*_1\cap F^*_2$ we define
\begin{align}\label{stochint}
\delta^*(h):=\int h(\eta^*,s,x)\eta^*(d(s,x))-\int h(\eta^*,s,x)\lambda(d(s,x)).
\end{align}
In particular,
\begin{align}\label{def2}
\delta^*(\I_{A\times(s,t]\times B}\I_{\bN^*\times C_n})
=\I_A(\eta^*)(\eta^*(((s,t]\times B)\cap C_n)-\lambda(((s,t]\times B)\cap C_n),
\end{align}
where $s<t$, $A\in \mathcal{N}_s$, $n\in\N$, and $B\in\mathcal{X}$.
Let $h\in F^*_1\cap F^*_2$ and define 
$\tilde{h}:\bN\times\BY\rightarrow\R$ by $\tilde{h}:=h$
on $\bN^*\times\BY$ and $\tilde{h}:=0$, otherwise. By Proposition \ref{l1},
$\tilde{h}$ is predictable. Since
$\BP(\eta\in\bN^*)=1$ we obtain from \eqref{tequal}
that $\delta^*(h)=\delta(\tilde h)$ $\BP$-a.s.
Therefore \eqref{l2i3} implies the isometry relation
\begin{align}\label{iso}
\BE \delta^*(h)^2=\BE \int h(\eta^*,s,x)^2\lambda(d(s,x))
\end{align}
for any $h\in F^*_1\cap F^*_2$. Since $F^*_1\cap F^*_2$ is dense
in $F^*_2$ we can extend $\delta^*$ 
to a linear operator from $F^*_2$ to $L^2(\BP)$.
Equation \eqref{iso} remains valid for arbitrary $h\in F^*_2$.

We now prove that $\delta$ extends the stochastic integral $\delta^*$.
Special cases of this result can be found
in \cite{Kab75} and \cite{NuViv90}. 
For $h:\bN \times\BY\rightarrow\R$,
the function $h^*:\bN^*\times\BY\rightarrow\R$ 
denotes the restriction of $h$ to $\bN^*\times\BY$.

\begin{theorem}\label{text2} 
Let $h\in L^2(\BP_\eta\otimes\lambda)$ such 
that $h^*$ is $\mathcal{P}$-measurable.
Then $\delta(h)=\delta^*(h^*)$ $\BP$-a.s.
\end{theorem}
{\sc Proof:} Since $\BP(\eta\in\bN^*)=1$, we have from Proposition \ref{text}
that $\delta(h)$ is defined.
By \eqref{tequal} and \eqref{stochint}
(and Proposition \ref{l1}) the assertion holds for any
$h\in F^*_1\cap F^*_2$. 
In the general case we may choose  
$h_k\in F^*_1\cap F^*_2$, 
$k\in\N$, such that $h_k\to h$ as $k\to\infty$ in $L^2(\BP_\eta\otimes\lambda)$.
Then $\delta^*(h^*_k)=\delta(h_k)$ converges to $\delta^*(h^*)$ 
in $L^2(\BP)$. Since $\delta$ is closed, this yields the assertion.\qed

\begin{remark}\label{rmnot}\rm 
Proposition \ref{l1} justifies our terminology for measurable functions
$h$ on $\bN\times\BY$ satisfying \eqref{predict}.
By this proposition, if $h$ is predictable then $h^*$ is 
$\mathcal{P}$-measurable.
Conversely, if $h^*$ is $\mathcal{P}$-measurable then there exists predictable
$\tilde{h}$ with $\tilde{h}^* = h^*$. 
If $h\in L^2(\BP_\eta\otimes\lambda)$ is predictable then
our notation $\int hd\hat\eta:=\delta(h)$ is justified by Theorem \ref{text2}. 
\end{remark}

\begin{remark}\label{rmususal}\rm A standard assumption in the
stochastic analysis literature is completeness of the
underlying filtration. Quite often one can find no
further comment on this technical (and sometimes annoying) hypothesis. 
In this paper we do not make this completeness assumption,
which is rather alien to point process theory.
\end{remark}

\section{Independent random measures}\label{sl2}
\setcounter{equation}{0}

Let $(\BY',\mathcal{Y}')$ be a Borel space
and $\beta$ be a $\sigma$-finite measure and diffuse measure on
$\BY'$.
Let $\mathcal{Y}'_0$ denote the system of all sets $B\in\mathcal{Y}'$
such that $\beta(B)<\infty$. In this section we consider
an {\em independent random measure} on $\BY'$
(see \cite{Ito56}) with {\em variance measure} $\beta$. This
is a family $\zeta':=\{\zeta'(B):B\in\mathcal{Y}'_0\}$
with the following three properties. First,
$\BE\zeta'(B)=0$ and $\BE\zeta'(B)^2=\beta(B)$ for any
$B\in\mathcal{Y}'_0$. Second, 
if $B_1,B_2,\ldots\in \mathcal{Y}'_0$ are pairwise
disjoint, then $\zeta'(B_1),\zeta'(B_2),\ldots$ are independent.
Third, if $B_1,B_2,\ldots\in \mathcal{Y}'_0$ are pairwise
disjoint and $B:=\cup B_n\in \mathcal{Y}'_0$ then
$\zeta'(B)=\sum_n\zeta'(B_n)$ in $L^2(\BP)$. By 
\cite[Theorem 4.1]{Kallenberg} the series also converges
almost surely.
Since $\beta$ is diffuse, it follows that the distribution of
$\zeta'(B)$ is infinitely divisible for any $B\in\mathcal{Y}'_0$,
see \cite[p.\ 81]{Kingman} for a closely related argument.
The L\'evy-Khinchin representation (see \cite[Corollary 15.8]{Kallenberg})
implies that
\begin{align}\label{LCR}
\log \BE e^{iu\zeta'(B)}=-a_Bu^2+
\int (e^{iuz}-1-iuz)\lambda(B,dz),
\quad u\in\R,
\end{align}
where $a_B\in\R$ and
$\lambda(B,\cdot)$ is a measure on $\R^*:=\R\setminus\{0\}$
satisfying $\int z^2 \lambda(B,dz)=\beta(B)$.
The measure $\lambda(B,\cdot)$ is the {\em L\'evy measure} of
$\zeta'(B)$ and is unique. 
We assume that $a_B=0$, so that $\zeta$ has no Gaussian
component. If $B\in \mathcal{Y}'_0$ is the disjoint union
of measurable sets $B_n$, $n\in\N$, then the independence
of the $\zeta'(B_n)$ and the uniqueness of the 
L\'evy measure implies that 
$\lambda(B,\cdot)=\sum^\infty_{n=1}\lambda(B_n,\cdot)$.
By a well-known result from measure theory (see \cite[p.\ 82]{Kingman})
there is a unique measure $\lambda$ on $\BY'\times\R^*$
such that $\lambda(B\times C)=\lambda(B,C)$ for all $B\in \mathcal{Y}'_0$ 
and all measurable $C\subset \R^*$. Hence equation \eqref{LCR}
can be rewritten as
\begin{align}\label{LCR2}
\log \BE e^{iu\zeta'(B)}=\int \I_B(x)(e^{iuz}-1-iuz)\lambda(d(x,z)),
\quad u\in\R,
\end{align}
whenever $\beta(B)<\infty$. 
By definition, 
\begin{align}\label{betla}
\int z^2\I\{x\in\cdot\}\lambda(d(x,z))=\beta(\cdot).
\end{align}
In particular, $\lambda$ is $\sigma$-finite.

Let us now consider a Poisson process $\eta$ 
on $\BY:=\BY'\times \R^*$ with intensity measure $\lambda$.
For any $B\in \mathcal{Y}'_0$ we define the Wiener-It\^o integral
\begin{align}\label{ze}
\zeta(B):=\int z\I_B(y)\hat\eta(d(y,z)).
\end{align}
Then $\zeta:=\{\zeta(B):B\in\mathcal{Y}'_0\}$
is an independent random measure with variance measure $\beta$.
We might think of a point of $\eta$ as being a point
in $\BY'$ with the second coordinate representing
its weight. Then the integral \eqref{ze} is the weighted
sum of all points lying in $B$, suitably compensated.
It follows from \eqref{LCR2} and basic properties of
$\eta$ (cf.\  \cite[Lemma 12.2]{Kallenberg} or \cite[Section 3.2]{Kingman})
that $\zeta(B)$ and $\zeta'(B)$ have the same
distribution for any $B\in\mathcal{Y}'_0$.
Henceforth it is convenient to work with $\zeta$
and the Poisson process $\eta$.

We now assume that $<'$ is a partial ordering on $\BY'$
satisfying the assumptions listed in the introduction,
where in \eqref{diffuse} the measure $\lambda$
has to be replaced with $\beta$. (If $y\in[y]$
for all $y\in\BY'$ this is strengthening
the diffuseness assumption on $\beta$.)
Then we can define a binary relation $<$ on $\BY=\BY'\times\R^*$
by setting $(y,z)<(y',z')$ if $y<'y'$. This relation
also satisfies our assumptions, where \eqref{diffuse}
comes from \eqref{betla} and the assumption on $\beta$.
The measurability of $(\mu,y)\mapsto \mu_y$ can be proved
using a measurable disintegration 
$\mu(d(y,z))=K(\mu,y,dz)\mu^*(dy)$, where $K$ is a 
kernel from $\bN\times\BY'$ to $\R^*$ and
$\mu\mapsto \mu^*$ is a measurable mapping from
$\bN=\bN(\BY)$ to $\bN(\BY')$ such that $\mu(\cdot\times\R^*)$
and $\mu^*$ are equivalent measures for all $\mu\in\bN$.

The stochastic integral of a predictable function 
$h:\bN\times\BY'\rightarrow\R$ against $\zeta$ is defined by
\begin{align}\label{hii}
\int h(\eta,y)\zeta(dy):=\int zh(\eta,y)\hat\eta(d(y,z))
\end{align}
provided that
\begin{align}\label{hi}
\BE \int h(\eta,y)^2 \beta(dy)
=\BE\int  z^2h(\eta,y)^2\lambda(d(y,z))<\infty.
\end{align}
Let $\mathcal{M}^2_\zeta\subset L^2(\BP)$ be the space
of all square integrable random variables $X$ given by
\begin{align}\label{eM}
X=\int h(\eta,y)\zeta(dy),
\end{align}
where the predictable function $h$ satisfies \eqref{hi}.
It follows from Proposition \ref{pquadcha}
that $\mathcal{M}^2_\zeta$ is a closed linear space.
Hence any $Y\in L^2(\BP)$ can be uniquely written as $Y=X+X'$,
where $X\in \mathcal{M}^2_\zeta$  and $X'\in L^2(\BP)$ is orthogonal
to $\mathcal{M}^2_\zeta$.
Decompositions of this type were first considered
by Kunita and Watanabe \cite{KuWa67}. The following theorem
makes this decomposition more explicit. 
We use a stochastic
kernel $J(y,dz)$ from $\BY'$ to $\R^*$ such that
\begin{align}\label{disl}
z^2\lambda(d(y,z))=J(y,dz)\beta(dy).
\end{align}
Such a kernel exists by a standard disintegration result 
(cf.\ \cite[Theorem 6.3]{Kallenberg} for a special case).

\begin{theorem}\label{tmrtKW} Let $f\in L^2(\BP_\eta)$
and define a predictable 
$h_f:\bN\times\BY'\rightarrow\R$ by
\begin{align}\label{DiN}
h_f(\eta,y)=\BE\Big[\int z^{-1}D_{(y,z)}f(\eta)
J(y,dz)\Big|\eta_{y}\Big].
\end{align}
Then $h_f$ satisfies \eqref{hi}
and we have  $\BP$-a.s.\ that 
\begin{align}\label{mrt3}
f(\eta)=\BE f(\eta)+\int h_f(\eta,y)\zeta(dy)+X',
\end{align}
where $X'\in L^2(\BP)$ is orthogonal to $\mathcal{M}^2_\zeta$.
\end{theorem}
{\sc Proof:} By Fubini's theorem applied to kernels we have
\begin{align*}
\BE \int h_f(\eta,y)^2 \beta(dy)
=\int\BE\Big(\int\BE[z^{-1}D_{(y,z)})f(\eta)|\eta_y]
J(y,dz)\Big)^2\beta(dy).
\end{align*}
Applying Jensen's inequality to the stochastic kernel
$J(y,dz)$ and using \eqref{disl} and \eqref{1} gives \eqref{hi}.
We now define $X'\in L^2(\BP)$ by
\begin{align}\label{Ntt}
X':=\int (\BE[D_{(y,z)}f(\eta)|\eta_y]
-zh_f(\eta,y))\hat\eta(d(y,z)).
\end{align}
Theorem \ref{tmrt2} implies \eqref{mrt3}. It remains to show
that $X'$ is orthogonal to $\mathcal{M}^2_\zeta$.
To this end we consider a random variable $X$ as given in \eqref{eM}.
By Proposition \ref{pquadcha},
\begin{align}\label{899}
\BE X X'=
\BE\int zh(y)(\BE[D_{(y,z)}f(\eta)|\eta_y]-zh_f(\eta,y))\lambda(d(y,z)).
\end{align}
We have
\begin{align*}
\BE&\int z^2h(y)h_f(\eta,y)\lambda(d(y,z))=
\BE\int h(y)h_f(\eta,y)\beta(dy)\\
&=\BE\iint h(y)\BE[z^{-1}D_{(y,z)}f(\eta)|\eta_y]
J(s,x,dz)\beta(d(s,x))\\
&=\BE\int zh(y)\BE[D_{(y,z)}f(\eta)|\eta_y]
\lambda(d(y,z)).
\end{align*}
Hence \eqref{899} implies $\BE X X'=0$, as claimed.\qed

\vspace{0.3cm}
Di Nunno \cite{Nunno07} proved Theorem \ref{tmrtKW}
for special (``core'') functions $f$ 
(and allowing also for a Gaussian part of $\zeta$)
in case $\BY'=\R_+\times\BX$, with $<'$ given as in Section \ref{spred}.
In the case where $J(y,\cdot)=\delta_1$
for $\beta$-a.e.\ $y$ (that is that $\zeta$ has only atoms
of size $1$), \eqref{mrt3} reduces to the Clark-Ocone type
formula \eqref{mr}.

The following result characterizes the class of square-integrable
stochastic integrals against $\zeta$.

\begin{corollary}\label{c09}
Let $f\in L^2(\BP_\eta)$ such that $\BE f(\eta)=0$. Then 
$f(\eta)\in \mathcal{M}^2_\zeta$ 
if and only if there is some predictable
$h:\bN\times\BY'\rightarrow\R$ satisfying \eqref{hi} such that
\begin{align}\label{exact}
\BE[D_{(y,z)}f(\eta)|\eta_y]=zh(\eta,y)\quad
\lambda\text{-a.e.}\, (y,z),\,\BP\text{-a.s.}
\end{align}
\end{corollary}
{\sc Proof:} Assume that \eqref{exact} holds. Then $h=h_f$ and 
the random variable $X'$ defined by \eqref{Ntt}
vanishes almost surely. Therefore Theorem \ref{tmrtKW}
shows that $f(\eta)$ can be written as a stochastic integral
against $\zeta$.

Assume conversely that $f(\eta)\in \mathcal{M}^2_\zeta$
and consider the decomposition \eqref{mrt3}. 
Since the orthogonal projection onto $\mathcal{M}^2_\zeta$
is unique, it follows that $X'=0$ $\BP$-a.s.
By definition  \eqref{Ntt} this means that \eqref{exact}
holds with $h:=h_f$.\qed

\section{Minimal variance hedging}\label{smeanva}
\setcounter{equation}{0}

We consider a Poisson process $\eta$ on $\BY:=\R_+\times\BX\times\BX'$, 
where $(\BX,\mathcal{X})$ and $(\BX',\mathcal{X}')$ are  Borel spaces. 
The partial ordering on $\BY$ is defined by
$(s,x,z)<(s',x',z')$ if $s<s'$. As always, the intensity measure
$\lambda$ of $\eta$ is assumed to satisfy \eqref{diffuse}.
Our aim in this section is to extend the results of Section \ref{sl2}
for the case $\BY'=\R_+\times\BX$. We replace $\R^*$ by the
general space $\BX'$ and the independent
random measure $\zeta$ by a more general $L^2$-valued signed random measure.
The special structure of $\BY'$ (and $\BY$) allows for a financial
interpretation of our results.
We consider a point $(s,x,z)$ of $\eta$ as representing 
a financial event at time $s$ of (asset) type $x$ and with mark $z$.
We let $\kappa:\bN\times\BY\rightarrow\R$ be a predictable function
and interpret $\kappa(\eta,s,x,z)$ as the size of the event $(s,x,z)$.
We assume that
\begin{align}\label{Mk}
\bar\beta(\cdot)
:=\BE\int \kappa(\eta,s,x,z)^2\I\{(s,x)\in\cdot\}\lambda(d(s,x,z))
\end{align}
is a $\sigma$-finite measure. The system of all measurable
$B\subset \R_+\times\BX$ such that $\bar\beta(B)<\infty$
is denoted by $\mathcal{Y}'_0$. For any $B\in\mathcal{Y}'_0$
we define by
\begin{align}\label{zet}
\zeta(B):=\int \kappa(\eta,s,x,z)\I_B(s,x)\hat\eta(d(s,x,z))
\end{align}
a square integrable random variable having $\BE\zeta(B)=0$. 
The stochastic integral
of a predictable $h:\bN\times\R_+\times\BX\rightarrow\R$ 
(here $(s,x)<(s',x')$ if $s<s'$) against $\zeta$ is defined by
\begin{align}\label{stint}
\int h(\eta,s,x)\zeta(d(s,x))
:=\int h(\eta,s,x)\kappa(\eta,s,x,z)\hat\eta(d(s,x,z))
\end{align}
provided that
\begin{align}\label{hiii}
\BE\int h(\eta,s,x)^2\kappa(\eta,s,x,z)^2 \lambda(d(s,x,z))<\infty.
\end{align}
We denote by $\mathcal{A}$ the set of all such 
predictable functions $h$.

\begin{remark}\label{rmfinance}\rm Let $\mathcal{X}_0$
denote the system of all $B\in\mathcal{X}$ such that
$[0,t]\times B\in \mathcal{Y}'_0$ for all $t\ge 0$. For
$B\in\mathcal{X}_0$ we can define the 
square integrable martingale (see Remark \ref{rpmart})
$$
\zeta_t(B):=\int \kappa(\eta,s,x,z)\I_{[0,t]}(s)\I_B(x)\hat\eta(d(s,x,z)),
\quad t\in[0,\infty].
$$
We interpret $\zeta_t(B)$ as the
(discounted) price of the assets in $B$ at time $t$. Note that
$\zeta_t(\cdot)$ is a signed measure on $\mathcal{X}_0$
in a $L^2$-sense. An element $h\in\mathcal{A}$ can
be interpreted as admissable portfolio investing the amount 
$h(\eta,s,x)$ in asset $x$ at time $s$.
Accordingly, if the bond price is constant, and $V_0\in\R$ then
$$
V_t:=V_0+\int\I_{[0,t]}(s) h(\eta,s,x)\zeta(d(s,x)), \quad t\in[0,\infty],
$$
is the value process of the self-financing portfolio 
associated with $h$ and an initial value $V_0$. 
\end{remark}

Let $f\in L^2(\BP_\eta)$. We interpret $f(\eta)$ as
a {\em claim} to be hedged (or approximated)
by a random variable of the form
$\BE f(\eta)+\int h(\eta,s,x)\zeta(d(s,x))$ with $h\in\mathcal{A}$.
A {\em minimal variance hedge} of $f(\eta)$ is then
a portfolio $h_f\in\mathcal{A}$ satisfying
\begin{align}\label{minhedge}\notag
\BE\Big(f(\eta)-\BE &f(\eta)-\int h_f(\eta,s,x)\zeta(d(s,x))\Big)^2\\
&=\inf_{h\in\mathcal{A}}
\BE\Big(f(\eta)-\BE f(\eta)-\int h(\eta,s,x)\zeta(d(s,x))\Big)^2.
\end{align}

\begin{remark}\label{rmfinance2}\rm Problem \eqref{minhedge}
requires us to minimize the quadratic risk among all
self-financing portfolios with initial value $\BE f(\eta)$.
We might also be interested in minimizing
\begin{align}\label{minhedge2}
\BE\Big(f(\eta)-c-\int h(\eta,s,x)\zeta(d(s,x))\Big)^2.
\end{align}
in $c\in\R$ and $h\in\mathcal{A}$. However, if $h_f\in\mathcal{A}$
solves \eqref{minhedge} then the pair $(\BE f(\eta),h_f)$
minimizes \eqref{minhedge2}.
\end{remark}

To solve \eqref{minhedge} we need to generalize the disintegration
\eqref{disl}. A kernel $J$ from $\bN\times\R_+\times\BX$ to 
$\BX'$ is called predictable, if
$(\mu,s,x)\mapsto J(\mu,s,x,C)$ is predictable for all
$C\in \mathcal{X}'$. In the next proof and also later
we use the generalized inverse $a^\oplus$ of a real number $a$.
It is defined by  $a^\oplus:=a^{-1}$ if $a\ne 0$ and
$a^\oplus:=0$ if $a=0$.

\begin{lemma}\label{l09} There exists a predictable
stochastic kernel $J$ from $\bN\times\R_+\times\BX$ to 
$\BX'$ such that
\begin{align}\label{disJ}
\kappa(\eta,s,x,z)^2\lambda(d(s,x,z))=J(\eta,s,x,dz)\beta(d(s,x))
\quad \BP\text{-a.s.},
\end{align}
where the random measure $\beta$ on $\R_+\times\BX$
is defined by
\begin{align}\label{beta}
\beta(\cdot):=\int \I\{(s,x)\in\cdot\}\kappa(\eta,s,x,z)^2\lambda(d(s,x,z)).
\end{align}
\end{lemma} 
{\sc Proof:} Define a measure $\bar\lambda$ on $\BY$ by
$$
\bar\lambda(d(s,x,z)):=\bar\kappa(s,x,z)\lambda(d(s,x,z)),
$$
where $\bar\kappa(s,x,z):=\BE\kappa(\eta,s,x,z)^2$.
Because the measure $\bar\beta=\bar\lambda(\cdot\times\BX')$ 
(see \eqref{Mk}) is assumed $\sigma$-finite and $\BX'$ is Borel,
there is a stochastic kernel $\bar{J}$ from $\R_+\times\BX$ to $\BX'$
such that 
$$
\bar\lambda(d(s,x,z))=\bar{J}(s,x,dz)\bar\beta(d(s,x)).
$$
It follows that
\begin{align}\label{4.51}
\kappa(\eta,s,x,z)^2\lambda(d(s,x,z))
=\bar\kappa(s,x,z)^{\oplus}
\kappa(\eta,s,x,z)^2 \bar{J}(s,x,dz)\bar\beta(d(s,x))
\quad \BP\text{-a.s.}
\end{align}
In particular the random measure $\beta$ defined by \eqref{beta}
coincides a.s.\ with $g(\eta,s,x)\bar\beta(d(s,x))$,
where
$$
g(\mu,s,x):=\int\bar\kappa(s,x,z)^{\oplus}\kappa(\mu,s,x,z)^2 \bar{J}(s,x,dz).
$$
We now define
$$
J(\mu,s,x,dz):=g(\mu,s,x)^{-1}\bar\kappa(s,x,z)^{\oplus}\kappa(\mu,s,x,z)^2
\bar{J}(s,x,dz),
$$
if $g(\mu,s,x)>0$. Otherwise we let $J(\mu,s,x,\cdot)$ equal some fixed
probability measure on $\BX'$. Then $J$ is predictable
and \eqref{4.51} implies \eqref{disJ}.\qed

\vspace*{0.3cm}
As in Section \ref{sl2} we let $\mathcal{M}^2_\zeta$ 
denote the space of all square integrable random variables
that can be written as a stochastic integral \eqref{stint}.

\begin{theorem}\label{tminvar} Let $f\in L^2(\BP_\eta)$
and define
\begin{align}\label{DiN2}
h_f(\eta,s,x)=
\int\kappa(\eta,s,x,z)^\oplus\BE[D_{(s,x,z)}f(\eta)|\eta_{s-}]J(\eta,s,x,dz),
\end{align}
where the stochastic kernel $J$ is as in Lemma \ref{l09}.
Then $h_f\in\mathcal{A}$ and \eqref{minhedge} holds.
Moreover, we have for any $t\in[0,\infty]$ that $\BP$-a.s.\ 
\begin{align}\label{mrt32}
\BE[f(\eta)|\eta_t]=\BE f(\eta)
+\int\I_{[0,t]}(s) h_f(\eta,s,x)\zeta(d(s,x))+N_t,
\end{align}
where $(N_t)$ is a square integrable martingale such
that $N_\infty$ is orthogonal to $\mathcal{M}^2_\zeta$.
\end{theorem} 
{\sc Proof:} Clearly $h_f$ is predictable. The integrability
condition \eqref{hiii} can be checked exactly as in the
proof of Theorem \ref{tmrtKW}. 
We can now proceed as in the proof of Theorem \ref{tmrtKW}
to derive the representation 
\begin{align}\label{mrt38}
f(\eta)=\BE f(\eta)+\int h_f(\eta,s,x)\zeta(d(s,x))+X',
\end{align}
where $X'\in L^2(\BP)$ is orthogonal to $\mathcal{M}^2_\zeta$.
This orthogonality implies \eqref{minhedge}. 
Let $t\ge 0$ and define $N_t:=\BE[X'|\eta_t]$. 
Taking conditional
expectations in \eqref{mrt38} and using Remark \ref{rpmart}
yields \eqref{mrt32}.\qed

\vspace{0.3cm}
The next result characterizes the claims that
can be perfectly hedged. The proof is an obvious
generalization of the proof of Corollary \ref{c09}.

\begin{theorem}\label{texacth}
Let $f\in L^2(\BP_\eta)$. Then 
\eqref{minhedge} vanishes if and only if there is some $h\in\mathcal{A}$ such that
\begin{align}\label{exacthedge}
\BE[D_{(s,x,z)}f(\eta)|\eta_{s-}]=\kappa(\eta,s,x,z)h(\eta,s,x)\quad
\lambda\text{-a.e.}\, (s,x,z),\,\BP\text{-a.s.}
\end{align}
In this case we have
$h(\eta,s,x)=h_f(\eta,s,x)$ for $\bar\beta$-a.e.\ $(s,x)$ and $\BP$-a.s.
\end{theorem}

In the remainder of this section we assume that
$\BX=\N$, that is, we assume that there are only
countably many assets. For any $j\in\N$ we define
a measure $\lambda_j$ on $\R_+\times\BX'$ by
$$
\lambda_j:=\iint\I\{(s,z)\in\cdot\}\lambda(ds\times\{j\}\times dz).
$$
Because $\lambda$ is $\sigma$-finite all measures $\lambda_j$
must be $\sigma$-finite as well.
Hence there exist $\sigma$-finite kernels $J_j$ from
$\R_+$ to $\BX'$ and $\sigma$-finite measures 
$\mu_j$ on $\R_+$ satisfying
$$
\lambda_j(d(s,z))=J_j(s,dz)\mu_j(ds),\quad j\in\N.
$$
The predictable function $\kappa$ is assumed to satisfy
$$
\BE\int\kappa(\eta,s,j,z)^2\lambda_j(d(s,z))<\infty,\quad j\in\N.
$$
This implies the $\sigma$-finiteness of the measure \eqref{Mk}.
The kernel $J$ of Lemma \ref{l09} is given by
$$
J(\mu,s,j,dz)=\Big(\int\kappa(\mu,s,j,z)^2 J_j(s,dz)\Big)^{-1}
\kappa(\mu,s,j,z)^2 J_j(s,dz)
$$
whenever $\int\kappa(\mu,s,j,z)^2 J_j(s,dz)>0$. If $f\in L^2(\BP_\eta)$
then, according to Theorem \ref{tminvar}, 
the minimal variance hedge $h_f$ of $f(\eta)$ can be computed as 
\begin{align}\label{DiN7}
h_f(\eta,s,j)=\Big(\int\kappa(\eta,s,j,z)^2 J_j(s,dz)\Big)^{\oplus}
\int\kappa(\eta,s,j,z)\BE[D_{(s,j,z)}f(\eta)|\eta_{s-}]J_j(s,dz).
\end{align}

\begin{example}\label{exfin}\rm Assume that $\BX'=\R^*$ 
and that
\begin{align*}
\int z^2\lambda_j([0,t]\times dz)<\infty,
\quad t\in\R_+,j\in\N. 
\end{align*}
Assume further that $\kappa(\eta,s,j,z)=\kappa_j(\eta,s)z$, 
for some predictable processes $\kappa_j$, $j\in\N$.
For any $h\in\mathcal{A}$ we then have
$$
\int h(\eta,s,j)\zeta(d(s,j))=
\sum_{j\in\N}\int h(\eta,s,j)\kappa_j(\eta,s)d\zeta_j(s)
\quad \BP\text{-a.s.},
$$
where $\zeta_j(t):=\iint \{s\le t\}z\hat\eta(ds\times\{j\}\times dz)$, 
$t\ge 0$, are independent square integrable processes with 
independent increments and mean $0$ (and no fixed jumps).
Assume now moreover, that $\lambda_j(d(s,z))=ds\nu_j(dz)$
for measures $\nu_j$ on $\R^*$, so that the
$\zeta_j$ are square integrable L\'evy martingales.
Then we can choose $J_j(s,dz)=\nu_j(dz)$ and \eqref{DiN7} simplifies to 
\begin{align}\label{DiN9}
h_f(\eta,s,j)=\kappa(\eta,s,j)^{\oplus}\Big(\int z^2\nu_j(dz)\Big)^{-1}
\int z\BE[D_{(s,j,z)}f(\eta)|\eta_{s-}]\nu_j(dz).
\end{align}
This is the main result in \cite{BNLOk03}. 
In fact, the model in \cite{BNLOk03} allows the processes
$\zeta_j$ to have a Brownian component
but considers only finitely many non-zero measures $\nu_j$.
\end{example}

\vspace{0.3cm}
\noindent
{\bf Acknowledgements:} We wish to thank Yuri Kabanov for giving
several comments on an early draft of this paper.

\end{document}